\documentclass[10pt,openright]{article}

\usepackage{amsfonts}
\usepackage{amssymb}
\usepackage{amsmath}
\usepackage{enumerate}

\newtheorem{Th}{\bf THEOREM}[section]

\newtheorem{Def}{\bf DEFINITION}[section]

\newtheorem{Pro}{\bf PROPOSITION}[section]

\newtheorem{St}{\bf STATEMENT}[section]

\newtheorem{Rq}{\bf REMARK}[section]

\newfont{\grandsy}{cmsy10}
\def\SS{{\grandsy x}}

\title{Set-valued differentiation as an operator}

\author{SERGUEI N. SAMBORSKI}

\date{}

\begin{document}

\maketitle

\begin{center}
Department de Mathematics, University de Caen, 14032 CAEN cedex, France

e-mail : samborsk@math.unicaen.fr
\end{center}

\vskip 1,5cm

\begin{abstract}
We introduce real vector spaces composed of set-valued maps on an open set. They are also complete metric spaces, lattices, commutative rings.
The set of differentiable functions is a dense subset of these spaces and the classical gradient may be extended in these spaces as a closed
operator. If a function $f$ belongs to the domain of such extension, then $f$ is locally lipschitzian and the values of extended gradient coincide
with the values of Clarke's gradient. However, unlike Clarke's gradient, our generalized gradient is a linear operator.
\end{abstract}

{\bf Key words} : functional metric spaces, functional lattices, extension of diffe\-rentiation, Clarke's gradient. 

{\bf AMS Subject Classification} : 26A24,28A15,46E05,54C35.\\ \\

\setcounter{section}{-1}

\section{Introduction.}

The assignment of the subset $\partial_Cf(x)$ to a locally Lipschitz function $f$ at a point $x$ (non smooth or set-valued or Clarke's
gradient) makes it possible to generalize various results of the classical differential calculus that involve the value of the
derivatives at points [2]. However, many important properties that are true in the classical analysis such as the linearity, the
Leibnitz formula for the derivative of a product fail under this pointwise approach. Moreover, these properties can not be even
rigourously formulated. For example we need to define the structures of a linear space in sets of set-valued maps to make it possible
to speak about the linearity of the differentiation. \\

In this paper we construct some linear spaces whose elements are set-valued maps. These spaces are also complete metric spaces,
lattices and commutative rings. Functions that are differentiable in the classical sense form a dense subset in each of these spaces.
The classical differentiation defines a preclosed operator on this subset in the constructed metric spaces. The closure of this
operator is an operator whose domain consists of locally Lipschitz functions. Actually we construct two such spaces using the general
scheme outlined in Appendix. Let us focus on the second of these spaces in which algebraic properties of classical differentiation
hold. This space consists of Riemann integrable functions, more exactly of equivalence classes of such a functions under the relation of coincidence almost everywhere. Consider the simpliest case of real-valued Riemann integrable functions. Assign to such a function
$\varphi$ and to a point $x$ the interval
$$[\lim_{y\rightarrow x}\inf \varphi(y), \lim_{y\rightarrow x}\sup\varphi(y)],$$
where $y$ is restricted to the set of points of continuity of the function $\varphi.$ It is easy to see that this interval is the same
for all functions equivalent to the function $\varphi.$ Denote this interval by $f(x),$ where $f$ be the equivalence class of
$\varphi.$ Thus the set-valued function $x\rightarrow f(x)$ is defined. This function  is single-valued almost everywhere. For
exemple, let $f_1$ be the equivalence class of the function

$$x \longmapsto
\begin{cases}
1 &$\hbox{ if } \ x $>$ 0$ \\
-1&$\hbox{ if } \ x $<$ 0.$ \\
\end{cases}$$

Then $f_1(0) = [-1,1].$ If $f_2$ is the equivalence class of the function ``$x\rightarrow 0$ if $x\not=0"$ then $f(0)=\{0\}.$ The
algebraic operations in this space are naturally defined by representatives, for example $f_1+(-f_1)=f_2.$\\

Let  $\boldsymbol{\partial}$  be the mentioned closure of the classical operator of differentiation and $\Delta$ be its domain. We show that if
$f\in\Delta$ then $\boldsymbol{\partial}$ 
 $f(x)=\partial_Cf(x)$ at every point $x,$ the equality means the coincidence of intervals in the
considered simpliest case or the coincidence of convex subsets in $\mathbb{R}^n$ in the general case.\\

The operator $\boldsymbol{\partial}$ 
defined by this way has the following property : its domain $\Delta$ is a linear subspace and $\boldsymbol{\partial}$  is a linear operator, i.e.
$$\boldsymbol\partial(\lambda f+g)(x) = (\lambda\boldsymbol\partial f+\boldsymbol\partial g)(x) \eqno(1.1)$$
for every $x$ and every $f,g\in\Delta, \lambda\in\mathbb{R}.$ \\

The domain $\Delta$ is a commutative subring and a sublattice as well. The corres\-ponding equalities for products and lattice
operations are true. This allows to cover the most important applications of the set-valued differentiation to various variational
problems.\\

In view of the well-known ``nonlinearity" of Clarke's gradient let us explain for a classical example the difference between the
point wise Clarke's approach and the operator approach in this article. Let $x\in\mathbb{R}, f(x)=|x|, g(x)=-|x|.$ In the Clarke's
theory [2], as well as under our approach one has $\partial_cf=\boldsymbol\partial f=f_1, \partial_cg=\boldsymbol\partial g=-f_1$ where $f_1$ is the
set-valued function from the example regarded above. Thus $\boldsymbol\partial f+\boldsymbol\partial g=0$ and the equality (1.1) is true for every $x.$
According to Clarke's approach only the following inclusion holds  :
$$\partial_C(f+g)(x)\subset \partial_Cf(x) +\partial_Cg(x),$$
where the addition on the right-hand side is the addition of convex sets. This inclusion can be obtained as a simple corollary of the
equality (1.1).\\

Note that the Clarke's definition of the set-valued gradient is applicable to every locally lipschitzian function while the
differentiation constructed in this paper has a smaller domain.
However, in the classical calculus the deepest results have been obtained for differentiable functions with supplementary conditions
for their derivatives (such as the continuity of derivatives, etc.). It corresponds, often implicitly, to the consideration of the
differentiation as an operator, for example in Banach spaces of continuous functions. The author assumes that the similar situation
takes place for the set-valued differentiation.\\

\section{Preliminaries and Notations.}

We fix an euclidian scalar product $<\cdot,\cdot>$ and a norm $\parallel\cdot\parallel$ in $\mathbb{R}^n.$ We denote by $B(x,r)$ the
open ball in $\mathbb{R}^n$ of radius $r$ centred in $x.$ If $A\subset\mathbb{R}^n$ is a subset, then $A_\varepsilon$ is its
$\varepsilon$-neighborhood. The Hausdorff distance between closed bounded subsets $A$ and $B$ of $\mathbb{R}^n$ is the number
$${\rm Inf}\{\varepsilon\geq 0 \mid A\subset B_\varepsilon \hbox{ and } B\subset A_\varepsilon\}.$$
By $gr f$ we denote the graph of a map $f,$ if $f$ is a set-valued map with the domain $X,$ then $gr f = \cup\{(x,\xi)\mid x\in X,
\xi\in f(x)\}.$ By $cl A$ we denote the closure of a subset $A\subseteq\mathbb{R}^n$ and by $co A$ the convex hull of $A.$ Let $f$ be a
locally lipschitzian real function on an open subset $X\subseteq\mathbb{R}^n.$ Denote by ${\cal D}(f)$ the subset of points of
differentiability of $f.$ Then
$$\partial_C f(x) = \bigcap_{\varepsilon>0}co\{\partial f(x')\mid x'\in{\cal D}(f)\cap B(x,\varepsilon)\}$$
when $\partial$ denote the classical gradient.\\

Let $E_1$ and $E_2$ be two metric spaces and $A$ be a single-valued operator from $\Delta\subset E_1$ to $E_2.$ The operator $A$ is
said to be preclosed if $\lim f_i=\lim g_i, f_i\in\Delta, g_i\in\Delta, Af_i\rightarrow F, Ag_i\rightarrow G$ imply $F=G.$ Let $A$ be a preclosed
opertor. One can put $\Delta' =\{f\in E_1\mid\exists$ a sequence $(f_i)_i$ in $\Delta$ converging to $f$ in $E_1$ such that
$(Af_i)_i$ converges in $E_2\}$ and to extend $A$ on $\Delta'$ by the rule $Af=\displaystyle\lim_{i\rightarrow\infty} Af_i.$ This extension
is said to be {\sl the extension of $A$ from $\Delta$ by closure.} The result of this extension is a closed single-valued operator,
i.e. its graph is a closed subset of $E_1\times E_2.$\\

Let $f$ be a bounded real valued function defined on a dense subset $X'\subset X.$ We denote by
$$(x\rightarrow f(x)\mid x\in X'\subset X)^*$$
the upper semicontinuous (u.s.c.) hull of $f$ on $X$ i.e. the function $x\rightarrow\lim\sup f(y),$ where $y\rightarrow x\in X, y\in
X'.$ Respectively
$$(x\rightarrow f(x)\mid x\in X'\subset X)_*$$
is the lower semicontinuous ($\ell$.s.c.) hull of $f$ on $X.$ If $f$ is defined on $X$ then we use the simplified notation $(f)^*$ for the
u.s.c. hull, respectively $(f)_*$ for $\ell$.s.c. hull.

 \section{Spaces $C_{cm}$ and $C_{ae}$.}

Let $X$ be an open bounded subset of $\mathbb{R}^n.$ Remind that a subset $Y\subset X$ is {\sl comeager,} if $X\setminus Y$ is a countable union of nowhere dense subsets in $X.$

\begin{Def} We denote by $C_{cm}(X,\mathbb{R}^m)$ the set of equivalence classes of bounded maps from $X$ to $\mathbb{R}^m$ for any of
which the set of continuity points is comeager and the equivalence relation means the coincidence on some comeager subset of $X.$
\end{Def}

(In this definition and further ``cm" means {\cal comeager).}
\medskip

Let $m=1.$ Then  in the set $C_{cm}(X,\mathbb{R})$ there are following algebraic structures :

\begin{itemize}
\item $C_{cm}(X,\mathbb{R})$ is a lattice, $f\leq g$ if for some (hence for every) representatives $\varphi\in f, \psi\in g$ the inequality
$\varphi(x)\leq\psi(x)$ holds for $x$ from some comeager subset of $X$ ;
\smallskip
\item  $C_{cm}(X,\mathbb{R})$ is a commutative ring, $f.g$ is the equivalence class of $\varphi.\psi,$ where $\varphi\in f, \psi\in g.$ 
For $m\geq 1$ :
\smallskip
\item $C_{cm}(X,\mathbb{R}^m)$ is a vector space over $\mathbb{R} : \lambda f+g$ is the equivalence class of $\lambda\varphi+\psi,$ where
$\varphi\in f,\psi\in g, \lambda\in\mathbb{R}$ ;
\smallskip
\item $C_{cm}(X,\mathbb{R}^m)$ is a commutative module over the ring $C_{cm}(X,\mathbb{R}).$
\end{itemize}

Let $f\in C_{cm}(X,\mathbb{R})$ and $k\in\mathbb{N}.$ There exist the greatest $k$-lipschitzian function that minorizes $f$ and the least
$k$-lipschitzian function that majorizez $f$ : 
\begin{eqnarray*}
f_k^- &={\rm Sup}\{\varphi\in Lip_k(X,\mathbb{R})|\varphi\leq f\}\\
f_k^+ &={\rm Inf}\{\varphi\in Lip_k(X,\mathbb{R})|\varphi\geq f\}.\\
\end{eqnarray*}
\\

\begin{Def}  {\sl (of the metric $s$ in the case $m=1$).}
\medskip
Let $f,g\in C_{cm}(X,\mathbb{R}).$ We set
$$s(f,g) =\underset{k\in\mathbb{N}}{\rm Sup} \max\{h(f_k^-,g_k^-), h(f_k^+,g_k^+)\},$$
where $h(\varphi,\psi)$ is the Hausdorff distance between the closures of the graphs in $cl X\times\mathbb{R}$ of bounded continuous
functions $\varphi$ and $\psi.$
\end{Def}

\begin{Th} {\sl The set $C_{cm}(X,\mathbb{R})$ endowed with the metric $s$ is a complete metric space and a conditionally complete
lattice. the subset $C^o(X,\mathbb{R})$ of bounded continuous functions on $X$ is a dense subset of $C_{cm}(X,\mathbb{R})$ both in the sense of metric
spaces and in the sense of lattices. The convergence induced by the metric $s$ on $C^o(X,\mathbb{R})$ is uniform.
}
\end{Th}

{\sl Proof.} We want to use Theorem A1 from Appendix with $\rho=h$ and ${\cal C}_k=Lip_k(X,\mathbb{R}).$ All conditions of the
construction of $\widetilde{\cal C}$ in Appendix are satisfied as it follows from Assertions 3.1-3.3 of the paper [7].\\

It remains to show that the metric sipace $(\widetilde{Lip(X,\mathbb{R}}),\widetilde h)$ in Theorem A1 is isometric to the space
$(C_{cm}(X,\mathbb{R}),s).$\\

Remind that a single-valued map $\phi : X\rightarrow\mathbb{R}$ is called to be {\sl quasicontinuous} if for any $x\in X$ and for any
$\varepsilon>0$ there exists an open subset ${\cal U}\subset X$ such that $x\in cl{\cal U}$ and for every $y\in{\cal U}$ we have
$|\phi(x)-\phi(y)|<\varepsilon,$ [4].\\

\begin{Pro}  {\sl In every class $f\in C_{cm}(X,\mathbb{R})$ there is a unique lower semi-continuous quasicontinuous representative
$f^-$ and there is a unique upper semi continuous quasicontinuous representative $f^+.$ Those representatives may be caracterized by
the following relation
$$(f^-)^* = f^+ \hbox{ and } (f^+)_* = f^- \eqno(2.1)$$
}
\end{Pro}

This proposition is prooved in [8] (with an unessential difference in the terminology). Let $x$ be a continuity point of the map
$f^-$ from Proposition 2.1. At the same time $x$ is a continuity point of $f^+$ and $f^-(x) = f^+(x).$ Because $f^-$ is
quasicontinuous, there is an open subset ${\cal U}$ such that $x\in cl{\cal U}$ and the values of $f^-$ in ${\cal U}$ are
sufficiently close to $f^-(x).$ The map $f_k^-$ from Appendix coincides with the Yosida transform of $f^- :$
$$f_k^-(x) = \underset{y\in X}{\rm Inf} (f^-(y)+k^{-1}||x-y||).$$
Hence for a sufficiently large $k \  f_k^-(x')$ is closed to $f^-(x)$ in a some point $x'$ closed to $x.$ The same conclusion is
valid for $f^+$ and $f_k^+.$ It prooves that $h(f_n^-,f_n^+)\rightarrow 0$ as $n\rightarrow\infty$ and the identical map $Lip(X,\mathbb{R})\rightarrow Lip(X,\mathbb{R})$
may be extended up to an isometry of metric spaces $(C_{cm}(X,\mathbb{R}),s)$ and $(\widetilde{Lip(X,\mathbb{R}}),\widetilde h).$
\medskip
The property of the order completness is prooved in [7]. \hfill$\square$
\\

\begin{Rq}{\rm It follows from Theorem 2.1 that the completion of the set $C^o(X,\mathbb{R})$ of bounded continuous functions with respect
to the metric $s$ is simultaneously the order completion in the sense of Dedekind-Mac Neille [1]. Such two-fold completion with
respect to $\mathbb{Q}$ is  a principal property of $\mathbb{R}$ (the equivalence of Cantor and Dedekind constructions). However such property does
not take place for known function metric spaces (with respect to subsets of continuous functions).}
\end{Rq}

Now we change the category properties on the corresponding measure properties, ``ae" means further {\sl almost everywhere.}
\\

\begin{Def} Denote by $C_{ae}(X,\mathbb{R}^m)$ the set of equivalence classes of bounded maps from $X$ to $\mathbb{R}^m$ any of which is
continuous almost everywhere and the equivalence relation means the coincidence almost everywhere in $X.$
\end{Def} 

In the set $C_{ae}(X,\mathbb{R}^m)$ there are natural structures of vector space and of module over $C_{ae}(X,\mathbb{R}),$ in $C_{ae}(X,\mathbb{R})$ there
are structures of lattice and commutative ring.
\\

If $f\in C_{ae}(X,\mathbb{R})$ and $k\in \mathbb{N}$ then there are well defined
\begin{eqnarray*}
f_k^- &={\rm Sup}\{\varphi\in Lip_k(X,\mathbb{R})\mid \varphi\leq f\} \\
f_k^+ &={\rm Inf}\{\varphi\in Lip_k(X,\mathbb{R})\mid \varphi\geq f\},\\
\end{eqnarray*}
where, of course, $\leq$ means the inequality almost everywhere.
\\

\begin{Def} {\sl (of the metric $r$ in the case $m=1)$.} Let $f, g\in C_{ae}(X,\mathbb{R}).$ We set 
$$r(f,g) = {\rm Sup}_{k\in\mathbb{N}}\max\{\delta(f_k^-,g_k^-), \delta(f_k^+,g_k^+)\},$$
where $\delta(\varphi,\psi) = h(\varphi,\psi) + \displaystyle\int_X|\varphi(x)-\psi(x)|dx$ for bounded continuous functions
$\varphi$ and
$\psi, \ h$ is the same metric as in Definition 2.2.
\end{Def} 

Determine an inclusion that is very important in what follows.
\\

Let $\varphi$ be an almost everywhere continuous function from $X$ to $\mathbb{R}^m.$ Then $\varphi$ is continuous on a dense subset and
consequantly $\varphi$ is continuous at points of some comeager subset [5]. Passing to the equivalence classes we have the well
defined map
$$i : C_{ae}(X,\mathbb{R}^m) \rightarrow C_{cm}(X,\mathbb{R}^m).$$
It is easy to see that this map is a continuous inclusion of metric spaces and an embedding of lattices.
\\

\begin{Pro}  {\sl Proposition 2.1 remains true if one replaces the set $C_{cm}(X,\mathbb{R})$ by te set $C_{ae}(X,\mathbb{R})$ and adds the
following property :}
$$\int_X f^+(x)dx = \int_X f^-(x)dx.$$
{\rm (The proof is evident in veiw of the introduced inclusion $i).$}
\end{Pro} 

\begin{Th}  {\sl The set $C_{ae}(X,\mathbb{R})$ endowed with the metric $r$ is a complete metric space. The subset $C^o(X,\mathbb{R})$ of
bounded continuous function is a dense subset in $C_{ae}(X,\mathbb{R}).$
}
\end{Th} 

{\sl Proof :} Follows from Theorem A1 and Proposition 2.2 by analogy with the proof of Theorem 2.1.
\\

\begin{Def}  Let $f\in C_{cm}(X,\mathbb{R}^m), x\in X.$ We call {\sl value} of $f$ in $x$ and we note it by $f(x)$ the following
convex subset in $\mathbb{R}^m$
$$f(x) = \bigcap_{\varepsilon >0} C_o\{\varphi(y)\mid y\in B(x,\varepsilon)\cap C_\varphi\},$$
where $\varphi\in f$ and $\varepsilon>0$ $C_\varphi$ is the subset of points of continuity of the map $\varphi.$
\end{Def} 

It is easy to see that $f(x)$ does not depend on the choice of representative $\varphi$ in $f$ but only on $f.$
\\

If $m=1$ then $f(x)$ coincides with the interval $[f^-(x),f^+(x)]$ and also may be defined as the interval
$$[\lim_{k\rightarrow\infty}f_k^-(x), \lim_{k\rightarrow\infty}f_k^+(x)].$$
In what if follows we shall identify the elements of sets $C_{cm}(X,\mathbb{R})$ or $C_{ae}(X,\mathbb{R}^m)$ with set-valued maps
$$f : x\rightarrow f(x)$$
from $X$ to the set of compact convex subsets in $\mathbb{R}^m.$
\\

Let us introduce metrics in the case of arbitrary $m.$
\\

For every $f\in C_{cm}(X,\mathbb{R}^m), \varphi\in f$ and $\xi\in\mathbb{R}^m$ the real function $x\rightarrow <\varphi(x),\xi>$ defines the equivalence
class in $C_{cm}(X,\mathbb{R}),$ we denote it by $<f(\cdot),\xi>.$ The same notation is used for $f\in C_{ae}(X,\mathbb{R}^m).$
\\

\begin{Def} {\sl (of the metrics $s$ and $r$ in the case $m>1).$} Let $f,g\in C_{cm}(X,\mathbb{R}^m).$ We set
$$s(f,g) = \underset{||\xi||=1}{\rm Sup} s(<f(\cdot),\xi>,<g(\cdot),\xi>).$$
Let $f,g\in C_{ae}(X,\mathbb{R}^m).$ We set
$$r(f,g) = \underset{||\xi||=1}{\rm Sup}  r(<f(\cdot),\xi>,<g(\cdot),\xi>).$$
\end{Def} 

\begin{Pro}  {\sl The sets $C_{cm}(X,\mathbb{R}^m)$ and $C_{ae}(X,\mathbb{R}^m)$ endowed with the metric $s$ and respectively with the
metric $r$ are complete metric spaces.
 }
\end{Pro} 

\begin{Pro} {\sl Suppose that a sequence $(F_n)_n$ in $C_{cm}(X,\mathbb{R}^n)$ converges to $F\in C_{cm}(X,\mathbb{R}^n),$ a sequence
$(x_n)_n$ of points of $X$ converges to $x\in X,$ a sequence $(\eta_n)_n$ of vectors in $\mathbb{R}^n$ converges to $\eta\in\mathbb{R}^n$ and
$\eta_n\in F_n(x_n)$ for every $n.$ Then $\eta\in F(x).$ The same is true if $C_{cm}(X,\mathbb{R}^n)$ is replaced by $C_{ae}(X,\mathbb{R}^n).$
}
\end{Pro} 

{\sl Proof :} If $n=1$ then the result follows immediately from the inequalities
$$\lim_{n\rightarrow\infty}F_n^-(x_n)\geq F^-(x), \lim_{n\rightarrow\infty}\sup F_n^+(x_n)\leq F^+(x),$$
where $[F^-(x), F^+(x)] = F(x).$
\\

Let $n>1.$ By the definition of the metric in the space $C_{cm}(X,\mathbb{R}^n)$ (respectively in $C_{ae}(X,\mathbb{R}^n))$ the sequence
$(<F_n,\xi>)_n$ converges to $<F,\xi>$ in $C_{cm}(X,\mathbb{R}^n)$ (respectively in $C_{ae}(X,\mathbb{R}^n))$ for every $\xi\in\mathbb{R}^n, ||\xi||=1.$
Because the subsets $F_n(x_n), F(x)$ are convex the inclusions $<\eta_n,\xi>\in <F_n(x_n),\xi>~,<\eta,\xi>\in <F(x),\xi>$  imply the desired result. \hfill$\square$
\\

\section{Gradient in $C_{cm}(X,\mathbb{R})$.}

Let $X$ be an open subset of $\mathbb{R}^n.$ Denote by $C^{1,0}$ the set of all differentiable real-valued functions posseding a continuous
bounded gradient in $X.$ $C^{1,0}$ is a dense subset in the metric space $C_{cm}(X,\mathbb{R})$ and the classical gradient is well defined on
$C^{1,0}$ as an operator from $C_{cm}(X,\mathbb{R})$ to $C_{cm}(X,\mathbb{R}^n).$
\\

\begin{Th} 
\begin{enumerate}[i)]
\item  {\sl The operator $\partial$ from $C^{1,0}\subset C_{cm}(X,\mathbb{R})$ to $C_{cm}(X,\mathbb{R}^n)$ is preclosed.
\smallskip
 Let us denote by $\boldsymbol\partial$ the closure of $\partial$ and by $\Delta_{cm}$ its domain. Then every element $f$ from
$\Delta_{cm}$ is a locally lipschitzian function and for every $x\in X$ the subset $\boldsymbol\partial,f((x)\subset\mathbb{R}^n$ coincides with the
subset $\partial_cf(x)$-the value of Clarke's gradient of $f$ in $x.$
}
\medskip
\item  $\Delta_{cm}$ {\sl is a sublattice of the lattice $C_{cm}(X,\mathbb{R})$ and for every $f,g\in\Delta_{cm}.$  $\boldsymbol\partial (\min(f,g))$ is
the equivalence class of the function
$$x\rightarrow
\begin{cases}
\partial f(x) &\hbox{ if } \ f(x)\leq g(x) \\
\partial g(x) &\hbox{ if } \ f(x)>g(x), \\
\end{cases}$$
where $x$ belongs to the comeager subset of points of differentiability of $f$ and $g.$ The analoguous formula is valid for $\boldsymbol\partial
(\max(f,g)).$
}
\end{enumerate}
\end{Th}

{\sl Proof : } i) We argue by reductio ad absurdum. Suppose that there exist two bounded sequence of differentiable functions
$(f_i)_i$ and $(g_i)_i$ converging to $f\in C_{cm}(X,\mathbb{R}),$ the sequences of their gradients $(\partial f_i)_i$ and $(\partial g_i)_i$
converge to $F$ and respectively to $G$ in $C_{cm}(X,\mathbb{R}^n)$ and $F\not= G.$ The intersection of the subsets of points of continuity of
$F$ and of $G$ is comeager, then there exist an open subset ${\cal U}$ and a real $\alpha>0$ such that the inequality
$||\xi-\eta||>\alpha$ holds for every $x\in{\cal U}, y\in{\cal U}, \xi\in F(x), \eta\in G(y).$ Let $x_0\in C(F)\cap
C(G)\cap{\cal U}.$ Then for every $\varepsilon>0$ there exists an open subset ${\cal U}'\subset{\cal U}$ such that $\forall x\in
{\cal U}', \forall y\in {\cal U}', \forall\xi\in F(x), \forall \eta\in G(y),$ we have
$$||F(x_0)-\xi||<\varepsilon \hbox{ and } ||G(x_0)-\eta||<\varepsilon.$$
Choosing $\varepsilon$ sufficiently small we see then there exist $\lambda\in\mathbb{R}^n, ||\lambda||=1,$ reels $\alpha_1,\alpha_2,
\alpha_1<\alpha_2,$ such that for every $x\in{\cal U}', \xi\in F(x), \eta\in G(x),$ the inequalities
$$<\xi,\lambda>\ <\alpha_1 \hbox{ and } <\eta,\lambda> \ >\alpha_2$$
hold.
\\

From the definition of the distance in $C_{cm}(X,\mathbb{R}^n)$ we have that the sequence $(<\partial f_i,\lambda>)_i$ converges to
$<F,\lambda>$ and the sequence $(<\partial g_i,\lambda>)_i$ converges to $<G,\lambda>$ in the space $C_{cm}(X,\mathbb{R}).$ Let $\varepsilon$
be equal $(\alpha_2-\alpha_1)/4.$ From Proposition 2.3 of the paper [8] it follows that there exist an open subset ${\cal
U}_1\subset{\cal U}'$ and a natural number $N$ such that the inequality
$$|<\partial f_n(x),\lambda>-<\varepsilon,\lambda>|<\varepsilon$$
holds for every $x\in{\cal U}_1, n>N_1$ and $\xi\in F(x).$ (Remark, that in [7,8] the space $C_{cm}(X,\mathbb{R})$ is denoted by $S(X),$
other notations are identic); Applying the same proposition to the open subset ${\cal U}_1$ and to the same $\varepsilon$ we obtain
that there exist an open subset ${\cal U}_2\subset{\cal U}_1$ and a natural number $N_2$ such that the inequality
$$|<\partial g_n(x),\lambda> - <\eta,\lambda>|<\varepsilon$$
holds for every $x\in{\cal U}_2, n>N_2$ and $\eta\in G(x).$ Taking into account the choice of $\varepsilon$ we obtain that the
inequalities
\begin{eqnarray*}
&<\partial f_n(x),\lambda> <\alpha'(={3\alpha_1+\alpha_2\over 4} )\\
&<\partial g_n(x),\lambda> >\alpha''(={3\alpha_2+\alpha_1\over 4})>\alpha' \\ 
\end{eqnarray*}
hold for every $x\in{\cal U}_2$ and every $n>N =\max(N_1,N_2).$ Let $x_0\in C(f)\cap{\cal U}_2.$ For every $\delta>0$ and for
sufficiently large $n$ we have the inequalities
$$|f_n(x_0)-f(x_0)|<\delta \hbox{ and } |g_n(x_0)- (x_0)|<\delta.$$
The mean-value theorem and (3.1) imply the following inequalities
\begin{eqnarray*}
&f_n(x_0+t\lambda)<\alpha't+\delta \\
&g_n(x_0+t\lambda)>\alpha''t-\delta\\
\end{eqnarray*}

for every $n>N$ on the interval $\{x_0+t\lambda\mid t\geq 0\}\cap{\cal U}_2.$ Choosing $\delta$ sufficiently small, we see that the
sequence $(s(f_i,g_i))_i$ can not tend to zero as $i\rightarrow\infty.$ This contradiction proves the preclousness of $\partial.$ The
following statement is easy to proove, it also is a particular case of Proposition 8.1 from [3].
\\

Denote by $h_-(f,g)$ the Hausdorff distance between the subgrafs of two bounded below semi-continuous functions $f$ and $g.$
\\

\begin{St} Let $f : X\rightarrow\mathbb{R}$ be a bounded below semi-continued function and $x\in X$ be such that the subderivative
$\partial_-f(x)$ of $f$ in $x$ is not empty. Let $\varepsilon>0, a\in\partial_-f(x).$ Then there exists $\delta>0$ such that for
every differentiable function $g : X\rightarrow\mathbb{R}$ satisfying $h_-(f,g)<\delta,$ the following inequality
$$\max\{||x-x'||,|f(x)-g(x')|, ||a-\partial g(x')||\}<\varepsilon$$ 
holds for some $x'\in X.$
\end{St}

Let $f_i\rightarrow f$ in  $C_{cm}(X,\mathbb{R}), \partial f_i\rightarrow F$ in $C_{cm}(X,\mathbb{R}^n).$ Then for every closed ball $B\subset X$
there exists a constant $K$ such that $||\partial f_i(x)||<K$ for every $x\in B$ and $i\in\mathbb{N}.$ Due to Statement 3.1 it implies that for every
$x\in B$ with non empty subset $\partial_-f(x)$ and for every $a\in\partial_-f(x)$ the inequality $||a||<K$ holds. From [6] it
follows that $f$ is lipschitzian in $B.$
\\

\begin{St}  {\sl Let $f\in\Delta_{cm}, x^o$ be a point of continuity of $\boldsymbol\partial f.$ Then $x^o$ is a point of
differentiability of $f$ and the equality $\boldsymbol\partial f(x^o) = \partial f(x^o)$ holds.
}
\end{St}

{\sl Proof :} Let $(f_i)_i$ be a sequence of continuously differentiable functions converging to $f$ in $C_{cm}(X,\mathbb{R}).$ Suppose that
the sequence $(\partial f_i)_i$ converge to $\boldsymbol\partial f$ in $C_{cm}(X,\mathbb{R}^n).$ Hence $(f_i)_i$ converges  $f$ for the uniform
distance and $(\partial f_i)_i$ tends to $\boldsymbol{\partial f}$ for the Hausdorff distance between graphs of functions. Let $\xi\in\mathbb{R}^n,
||\xi||=1.$ Let us localize the study in a sufficiently small neighborhood of the point $x^o.$ Introduce in this neighborhood local
coordinates $(x_1,x_2)$ corresponding to the product structure $\{\lambda\xi\mid\lambda\in]\alpha,\beta[\}\times
V\subset\mathbb{R}\times\mathbb{R}^{n-1}=\mathbb{R}^n,$ where $V$ is some open subset of $\mathbb{R}^{n-1}.$ Denote $\boldsymbol\partial f$ by $F,$ $<F,\xi>_i^- \hbox{ by } F_i^-, <F,\xi>_i^+ \hbox{ by } F_i^+$
(see 2 for the definition of the projections $f\rightarrow f_i^{\pm}.$) The continuity of $F$ in $x^o$ and the nature of the convergence of
the sequence $(\partial f_i)_i$ to $F$ are such that there exist a decreasing sequence of reals $(\alpha_i)_i$ tending to zero, an
increasing sequence of reals $(\beta_i)_i$ tending to zero with the following properties
\begin{eqnarray*}
&F_i^+(\tau,x_2^o) + \alpha_i> <\partial f_j(\tau, x_2^o),\xi> \\
&F_i^-(\tau, x_2^o) + \beta_i< <\partial f_j(\tau, x_2^o),\xi> \\
\end{eqnarray*}
for every $\tau \in ]\alpha,\beta[$ and ever $j<i.$ Let us pose
\begin{eqnarray*}
\varphi_i^+(x_1) &= \int_{x_1^o}^{x_1}(F_i^+(\tau,x_0^o)+\alpha_i)d\tau + f(x_1^o,x_2^o) \\
\varphi_i^-(x_1) &= \int_{x_1^o}^{x_1} (F_i^-(\tau, x_2^o)+\beta_i)d\tau + f(x_1^o,x_2^o).\\
\end{eqnarray*}
The sequence $(F_i^+(\cdot,x_2^o))_i$ of real functions is bounded and decreasing, the same is valid for the sequence
$(F_i^+(\cdot,x_2^o)+\alpha_i)_i.$ Hence both sequences converge in the $L_1$-norm to $F^+(\cdot,x_2^o).$ It means that the sequence
$(\varphi_i^+(\cdot))_i$ converges in the uniform metric to some continuous function $\varphi^+.$ On the interval
$[x_1^o,x_1^o+\lambda[$ the sequence $(\varphi_i^+)_i$ is decreasing and $\varphi_i^+(\cdot)>f_i(\cdot,x_2^o).$ The continuity of
$F^+$ in the point $(x_1^o,x_2^o)$ implies that the function $\varphi^+$ has the right-side derivative in $x_1^o$ equal to
$<F^+(x_1^o,x_2^o),\xi>.$ We have also the inequality
$$\varphi_i^+(x_1) \geq f(x_1,x_2^o)$$
for $x_1\geq x_1^o,$ where $x_1$ is sufficiently close to $x_1^o.$ By the same way we obtain that the sequence
$(\varphi_i^-(\cdot))_i$ converges to some continuous function $\varphi^-.$ This function has right-side derivative in $x_1^o$ which
is equal to $<F^-(x_1^o,x_2^o),\xi>.$ The inequality
$$\varphi_i^-(x_1)\leq f(x_1,x_2^o)$$
is valid for $x_1\geq x_1^o.$ Then the equalities $\varphi^+(x_1^o) = \varphi^-(x_1^o) = f(x_1^o,x_2^o)$ and $F^-(x_1^o,x_2^o) =
F^+(x_1^o,x_2^o)$ imply the differentiability of the function $f$ in $x^o$ in the direction $\xi.$ Because its derivative is equal to
$<F(x^o),\xi>$ we obtain the differentiability of $f$ in $x^o$ with the gradient $\partial f(x^o)$ equals to $F(x^o)=\partial
f(x^o).$ The statement is prooved. \hfill$\square$
\\

Let $f\in\Delta_{cm}$ ; then for every $x\in X$ the inclusion $\partial_Cf(x)\subseteq\boldsymbol\partial f(x)$ follows from the defintions of
the Clarke's gradient and of the value in $x$ of $\boldsymbol\partial f\in C_{cm}(X,\mathbb{R}^n).$ Let us proove the inclusion $\boldsymbol\partial f(x)
\subseteq\partial_Cf(x).$
\\

Statement 3.2 implies the following inclusion
$$\boldsymbol\partial f(x)=\bigcap_{\varepsilon>0}co\{\boldsymbol\partial f(x')\mid x'\in C(\boldsymbol\partial f)\cap
B(x,\varepsilon)\}$$
$$\subset\bigcap_{\varepsilon>0}co\{\partial f(x')\mid x'\in{\cal D}(f)\cap B(x,\varepsilon)\}=\partial_Cf(x).$$
Assertion i) of the theorem is prooved. ii) Let $f,g\in\Delta_{cm}.$ Then $X = \Omega_1\cup\Omega_2\cup\Gamma,$ where
$\Omega_1=\{x\in X\mid f(x)<g(x)$ or $g(x)<f(x)\}, \Omega_2$ is the interiour of the subset $\{x\in X\mid f(x)=g(x)\}$ and $\Gamma$
is its boundary. It is evident that the map in Assertion ii) of the theorem coincide with $\boldsymbol\partial(\min(f,g))$ on the subset
$\Omega_1\cup\Omega_2.$ However the subset $\Gamma$ is nowhere dense, hence the values of the map on it may be neglected in the proof
that the map belongs to the space $C_{cm}(X,\mathbb{R}^n)$ and in the formula for $\boldsymbol\partial(\min(f,g)).$
\\

Proove now that $\min(f,g)\in\Delta_{cm}.$ Let $(f_n)_n, (g_n)_n$ be two sequences of $C^{1,0}$-functions converging to $f$ and
respectively to $g$ in the space $C_{cm}(X,\mathbb{R}).$ Suppose also that the sequences $(\partial f_n)_n$ and $(\partial g_n)_n$ converge in
the space $C_{cm}(X,\mathbb{R}^n).$ Let us localize our consideration. Suppose that $f,g,f_n,g_n$ are lipschitzien with the Lipschitz
constant $k.$ Let $(\varepsilon_n)_n$ be a sequence of positive reels tending to 0 and $(h_n)_n$ be a sequence of $k$-lipschitzien
$C^1$-functions with the following properties~:

\begin{enumerate}[i)]
\item  outside of the $\varepsilon_n$-neighborhood of the subset $\phi_n=\{x\mid f_n(x)=g_n(x)\}$ the function $h_n$ coincide
with the function $\min(f_n,g_n)$ ;

\item  inside of the $\varepsilon_n$-neighborhood of the subset $\phi_n$ the function $h_n$ satisfy the inequality
$$\min\{<\ f_n(\cdot),\xi>,<\partial g_n(\cdot),\xi>-k\varepsilon_n\leq <\partial h_n(\cdot),\xi> \leq \max\{<\partial
f_n(\cdot),\xi>,$$

$<\partial g_n(\cdot),\xi>\}+ k\varepsilon_n$ for every $\xi\in \mathbb{R}^n, ||\xi||=1.$
\end{enumerate}

Then the sequence $(h_n)_n$ tends to $\min(f,g)$ as $n\rightarrow\infty$ and simultunously the sequence $<\partial h_n(\cdot),\xi>$ tends
to some limit in the space $C_{cm}(X,\mathbb{R}).$ Hence the sequence $(\partial h_n)_n$ tends to a limit in the space $C_{cm}(X,\mathbb{R}^n).$ The
closeness of the operator $\boldsymbol\partial$ implies that the function $\min(f,g)$ belongs to the subset $\Delta_{cm}$ and that the formula in
the assertion ii) is true. \hfill$\square$
\\

\section{Gradient in $C_{ae}$.}

The classical gradient $\partial$ is also well defined on the set $C^{1,0}$ of continuously differentiable functions with bounded
gradient as an operator from $C_{ae}(X,\mathbb{R})$ to $C_{ae}(X,\mathbb{R}^n).$
\\

\begin{Th}
\begin{enumerate}[i)]
\item  {\sl The operator $\partial$ with the domain $C^{1,0}$ is preclosed from $C_{ae}(X,\mathbb{R})$ to $C_{ae}(X,\mathbb{R}^n).$ Denote by
$\boldsymbol\partial$ the closure of $\partial$ and by $\Delta$ its domain. Then $\Delta\subset\Delta_{cm},$ therefore every $f\in\Delta$ is a
locally lipschitzian function and for every $x\in X$ there is the equality of convex subsets
$$\boldsymbol\partial f(x) = \partial_Cf(x).$$
}
\item  $\Delta$ {\sl is a linear subspace of the linear space $C_{ae}(X,\mathbb{R})$ and $\boldsymbol\partial$ is a linear operator : for every
$f,g\in\Delta,\lambda\in\mathbb{R}$ there is the equality
}
$$\boldsymbol\partial(\lambda f+g) = \lambda\boldsymbol\partial f+\boldsymbol\partial g.$$
\item $\Delta$ {\sl is a commutative subring of the ring $C_{ae}(X,\mathbb{R})$ and for every $f,g\in\Delta$ there is the equality
$$\boldsymbol\partial(f.g) = f.\boldsymbol\partial g+g.\boldsymbol\partial f.$$
}
\end{enumerate}
\end{Th}

{\sl Proof.  :}
\begin{enumerate}[i)]
\item Remind that there is the continuous inclusion $i : C_{ae}(X,\mathbb{R}^n) \rightarrow C_{cm}(X,\mathbb{R}^n)$ (\SS2), and the convergence
of a sequence $(f_n)_n$ to $f$ in the space $C_{ae}$ imply the corresponding convergence in the space $C_{cm}.$ Hence the precloseness of
$\partial$ in the space $C_{cm}$ implies the precloseness of $\partial$ in the space $C_{ae}.$ The rest follows from Assertion i) of
Theorem 3.1.
\item  Let $f,g$ belong to $\Delta.$ It means that there are sequences $(f_n)_n, (g_n)_n$ of functions from $C^{1,0}$ such that
$f_n\rightarrow f, g_n\rightarrow g, \partial f_n\rightarrow\boldsymbol\partial f, \partial g_n\rightarrow\boldsymbol\partial g.$ We want to proove that $f+g\in\Delta.$ However, we
can not use the sequence $(f_n+g_n)_n$ as an approximation of $f+g.$ It may be happen that the sequence  $(\partial f_n+\partial
g_n)_n$ does not converge in $C_{ae}(X,\mathbb{R}^n).$ Then we shall modify the sequences $(f_n)_n$ and $(g_n)_n$ in order to have a
``continuity of addition" for the images of the modified sequences under the map $\boldsymbol\partial.$
\end{enumerate}

{\bf LEMMA 4.1.} {\sl Let $f,g$ be two elements of the space $C_{cm}(X,\mathbb{R}).$ Then the sequence $(f_n^-+g_n^+)_n$ converge to $f+g$ in
$C_{cm}(X,\mathbb{R}).$
}
\\

{\sl Proof :} Let $x\in C_f\cap C_g.$ Then for every $\varepsilon>0$ there exists an open subset ${\cal U}$ such that $x\in{\cal U}$
and
$$|f(y)-f(x)| <\varepsilon, |g(y)-g(x)|<\varepsilon$$
for every $y\in{\cal U}.$ The functions $f$ and $g$ are bounded, hence there exists an open subset $V\subset{\cal U}$ and $N\in\mathbb{N}$
such that
$$|(f_n^- + g_n^+)(y) - (f+g)(x)|< 2\varepsilon $$ 
for every $y\in V$ and every $n>N.$
\\

Suppose now that $x\notin C_f, f(x) = [a_-,a_+]$ and $g(x)=[b_-,b_+] (b_-=b_+$ if $x\in C_g).$
\\
For every $\varepsilon>0$ and for every open subset $W$ containing the point $x$ there are two points $w_-$ and $w_+$ from the
subset $C_f\cap C_g\cap W$ such that the inequality
$$(f+g)(w_-)-\varepsilon <a_-+b_+ < (f+g)(w_+)+\varepsilon$$
holds.
\\

For a sufficiently large $n\in\mathbb{N}$ and for a sufficiently small open subset $W$ the values of $f_n^-$ in $W$ are sufficiently closed
to $a_-$ and the values of $g_n^+$ in $W$ are sufficiently closed to $b_+.$ Hence we have the inequality
$$(f+g)(w_-)-2\varepsilon < f_i^-(x) + g_i^+(x) < (f+g)(w_+)+2\varepsilon.$$
This inequality together with (4.1) show that all limit values of the sequence $f_n^-+g_n^+$ are determined by the convex hulls of
limits of values of $f+g$ on the subset $C_+\cap C_g.$ From the definition of addition in $C_{cm}(X,\mathbb{R})$ we obtain that the Hausdorff
distance between the graphs of $f_n^-+g_n^+$ and the graph of $f+g$ (as a set valued map) tends to 0 as $n\rightarrow\infty.$ To conclude
the proof of the lemma it remains to refer to the following statement.
\\

\begin{St}{\sl Let $f\in C_{cm}(X,\mathbb{R}), (f_n)_n$ be a sequence of continuous real functions on $X.$ Suppose that the
Hausdorff distance between the graphs of $f_n$ and $f$ tends to 0 as $n\rightarrow\infty.$ Then $(f_n)_n$ converges to $f$ for the distance
of the space $C_{cm}(X,\mathbb{R}).$
}
\end{St}

{\sl Proof :} Assume the converse. Then the sequence $(f_n)_n$ does not converge in the space $C_{cm}(X',\mathbb{R})$ where $X'$ is a
bounded open subset of $X.$ Hence there exist a subsequence $(j_n)_n$ of natural numbers, a point $x\in cl X,$ a positive real number
$\delta>0$ and a real $b$ such that in the open ball $B(x,\delta)$ we have either $(f_n)_{j_n}^-<b-\delta$ and simultaneously
$f^+(x)>b+\delta$ or $(f_n)_{j_n}^+>b+\delta$ and simultanously $f^-(x)<b-\delta$ or these two cases at the same time.
\\

We consider the first case, the second case may be treated similarly and the third one is included in one of the preceding case. It
is easy to see that $j_n\rightarrow\infty$ as $n\rightarrow\infty.$ The function $f^+$ is quasi-continuous hence there exists an open subset
$V\subset B(x,\delta)$ such that all values $f^+(w)$ are suffciently closed to $f^+(x)$ for $w\in V.$ For sufficiently large $n$ this
contradicts to the convergence of Hausdorff distances declared in the statement.
\hfill$\square$ \\

Continue to proove Assertion ii) of the theorem. Assume, by contradiction, that there are $f\in\Delta, g\in\Delta$ such that
$f+g\notin\Delta.$ In view of Assertion i) it means that there exist an open subset ${\cal U}\subset X, \xi\in\mathbb{R}^n, ||\xi||=1$ such
that for every sequence $(\varphi_n)_n$ of functions from $C^{1,0}$ converging to $f+g$ the sequence $(<\partial\varphi_n,\xi>)_n$
does not converge in the space $C_{ae}({\cal U},\mathbb{R}).$ Let $F = <\partial f,\xi>, G = <\partial g,\xi>.$ 
\\

Let us replace the Lipschitz
functions $F_n^-$ and $G_n,^+$ by some smooth functions $F_n$ and $G_n$ suficiently closed to $F_n^-$ and $G_n^+$ in the uniform
metric. We choose this closeness in order to provide the convergence of $(F_n+G_n)_n$ to $F+G$ (see Lemma 4.1) and to conserve the
monotonicity of the sequences $(F_n)_n, (G_n)_n.$ 
\\

Let us localize our study in a sufficiently small open subset ${\cal U}\subset X.$ Introduce in ${\cal U}$ a local coordinates
$(x_1,x_2)$ corresponding to the product structure
$${\cal U} = \{\lambda\xi\mid\lambda\in]\alpha,\beta[\subset\mathbb{R}\}\times V\subset\mathbb{R}\times\mathbb{R}^{n-1}=\mathbb{R}^n$$
where $\alpha,\beta\in\mathbb{R}, V$ is some open subset of $\mathbb{R}^{n-1}.$ Let
\begin{eqnarray*}
&\varphi_i(x_1,x_2) = \int_\gamma^{x_1} F_i(\tau,x_2) d\tau+f_i(\gamma,x_2)\\
&\psi_i(x_1,x_2) = \int_\gamma^{x_1} G_i(\tau,x_2)d\tau + g_i(\gamma,x_2) \\
\end{eqnarray*}
where $\gamma\in]\alpha,\beta[, (\gamma,x_2)\in C_f\cap C_g.$
\\

Then we have the equalities
$$<\partial\varphi_i,\xi> = F_i, <\partial\psi_i,\xi> = G_i.$$
The sequence of continuous functions $(F_i)_i$ is increasing and bounded, hence it has a limit in the space $L_1({\cal U}),$ this
limit is just the function $F^-.$ The decreasing bounded sequence of continuous functions $(G_i)_i$ has a limit in the space
$L_1({\cal U}),$ it is just the function $G^+$ (remind that according to Proposition 2.1 $F = (F^-,F^+)$ and $G = (G^-,G^+)).$ In
particular, for every $\varepsilon>0$ there exists a $N\in\mathbb{N}$ such that the inequalities
$$\int_\gamma^{x_1} |F_i-F_j|(\tau,x_2)d\tau<\varepsilon, \int_\gamma^{x_1}|G_i-G_j|(
\tau,x_2)d\tau<\varepsilon$$
hold for every fixed $(x_1,x_2)$ and for every $i,j>N.$ It means that the sequences $(\varphi_i)_i, (\psi_i)_i$ 	are Cauchy sequences
for the uniform metric in ${\cal U}.$ Hence they converge. Denote their limits by $\varphi^-$ and $\psi^+.$ We have that
$\varphi^-\leq f$ and $\psi^+\geq g.$ By changing the roles of $f$ and $g$ we obtain two limit functions $\varphi^+\geq f$ and
$\psi^-\leq g.$ However $F$ and $G$ are continuous almost everywhere and $F^- = F^+, G^-=G^+$ almost everywhere in ${\cal U}.$ Hence
$\varphi^-=\varphi^+=f$ and $\psi^-=\psi^+=g.$
\\

Finally, the sequence $(\varphi_i+\psi_i)_i$ converges uniformly to $f+g$ and the sequence of terms
$$<\partial(\varphi_i+\psi_i),\xi> = <\partial\varphi_i,\xi> + <\partial\psi_i,\xi> = F_i+G_i$$
converges to $F+G$ in $C_{cm}({\cal U},\mathbb{R})$ by Lemma 4.1 and simultaneously in $L_1({\cal U}).$ It follows that the sequence
$(F_i+G_i)_i$ converge to $F+G$ in $C({\cal U},\mathbb{R}).$ This contradiction proves the assertion ii) of the theorem.\\

iii) The following particular case of the chain rule is easy to proove using approximative sequences from $C^{1,0}$ and the closeness
of the operator $\boldsymbol\partial.$
\\

\begin{St}  Let I be an interval in $\mathbb{R}, \phi : I\rightarrow\mathbb{R}$ be continuously differentiable function with bounded derivative on
I. Let $f\in\Delta, \rm Im f\subset I.$ Then  the composite function $\phi\circ f$ belongs to $\Delta$ and the following equality is
true :
$$\boldsymbol\partial(\phi\circ f) = \phi'(f(\cdot))\boldsymbol\partial f.$$
\end{St}
{\sl Proof.}
Let $f,g\in\Delta, \alpha>0.$ There are real constantes $C_1, C_2$ such that the values of the functions $f+C_1, g+C_2$ belong to a
bounded interval $[\alpha,\beta]$ with some $\beta\in\mathbb{R}.$ Apply the statement 4.2 with $\phi = \log$ to the product
$(f+C_1).(g+C_2).$ Using the assertion ii) of the theorem we see that $\log(f+C_1)(g+C_2)$ belongs to $\Delta,$ hence the product
$fg$ also belongs to $\Delta.$ The announced formula follows from the properties of the logarithmic function.

The theorem is prooved. \hfill$\square$
\\

Such results as Extremum Conditions, Mean Value Theorem and others are valid for the functions differentiable in the sens of this
paper. This follows from Theorems 3.1 and 4.1 because these results hold for the Clarke's gradient. However,  they can be obtained
immediately as consequences of well known theorems of the classical diffferential calculs by the passage to limit and by using of
Proposition 2.3.
\\

Suppose,  for example, that $f\in\Delta$ has a local extremum in a point $x_0\in X.$ Let $(f_n)_n$ be a sequence of differentiable
functions converging to $f$ in $C_{ae}(X,\mathbb{R})$ such that the sequence $(\boldsymbol\partial f_n)_n$ converges to $\partial f$ in
$C_{ae}(X,\mathbb{R}^n).$ Every function $f_n$ has a local extremum in some $x_n\in X,$ hence $\partial f_n(x_n)=0.$ Choosing a subsequence
$(x_{n_k})_k$ converging to $x\in X$ and applying Proposition 2.4 we obtain that $0\in \boldsymbol\partial f(x_0).$
\\

The same reasoning is valid for demonstrations of Mean Value Theorem and for other similar results.
\\

We complete this paper by a generalisation of classical results concerning the passage to limit under the sign of differentiation.
\\

Remind that if $F, G$ are two bounded set-valued maps then $h(F,G)$ denotes the Hausdorff distance between the closures of the graphs
of $F$ and $G$ in $cl X\times\mathbb{R}^n.$
\\

\begin{Pro}{\sl Let $(f_n)_n$ be a sequence of functions from $\Delta, x_0$ be a point of $X$ such that the sequence of
reals numbers $(f_n(x_0))_n$ converges in $\mathbb{R}.$ Suppose that the sequence $(\boldsymbol\partial f_n )_n$ converges in the metric $h$ to $F\in
C_{ae}(X,\mathbb{R}^n).$ Then the sequence $(f_n)_n$ converges uniformly to a function $f\in\Delta$ and $\boldsymbol\partial f=F.$
}
\end{Pro}

The proof makes use to the same standard local arguments as in the proof of assertion ii) of Theorem 4.1, Statement 4.1, the
closeness of the operator $\boldsymbol\partial$ and the following property : if $(\varphi_n)_n$ be a sequence of real functions from $\Delta$
converging to $\varphi$ in the metric $h$ and if $\varphi\in C_{ae}(X,\mathbb{R}),$ then $\varphi_n$ converges to $\varphi$ in $L_1(X,\mathbb{R}).$
\\

{\bf Appendix. Inductive limits of metric lattices.}
\\

Let ${\cal C}$ be a lattice with partial order $\leq$ endowed with a metric $\rho.$ We assume that the following conditions hold :\\

\begin{enumerate}[i)]
\item $\forall a, b\in{\cal C}, a<b \ \ \exists c\in {\cal C}$ such that $a<c<b.$

\item $\forall a,b,c\in{\cal C}$ such that $a<c<b, \rho(a,c)\leq\rho(a,b)$ and $\rho(c,b)\leq\rho(a,b).$

\item any bounded monotone sequence of elements of ${\cal C}$ is a Cauchy sequence in the metric $\rho.$
\end{enumerate}

Let $({\cal C})_{i\in\mathbb{N}}$ be a sequence of subsets of ${\cal C}$ such that ${\cal C} = \displaystyle\bigcup_{i\in\mathbb{N}}{\cal C}_i$ and
${\cal C}_n\subseteq{\cal C}_{n+1}$ for any $n\in\mathbb{N}.$
\\

Suppose that $\forall n\in\mathbb{N}$ the subset ${\cal C}_n$ with the partial order and the metric induced from ${\cal C}$ is a conditionally
complete lattice and a complete metric space simultaneously. The conditional completeness of ${\cal C}_n$ makes it possible to define
projectors from ${\cal C}$ to ${\cal C}_n$ by the rules

\begin{eqnarray*}
f\rightarrow f_n^- &={\rm Sup}\{\varphi\in{\cal C}_n\mid\varphi\leq f\} \\
f\rightarrow f_n^+ &={\rm Inf}\{\varphi\in{\cal C}_n\mid\varphi\geq f\}.\\
\end{eqnarray*}
The last condition on $({\cal C},\leq,\rho)$ is formultated in terms of these projectors :

{\sl for every $n\in\mathbb{N},$ the mappings $f\rightarrow f_n^-, f\rightarrow f_n^+$ are uniformly continuous.
}
\\

\begin{Def}We define $\widetilde{\cal C}$ as the set of pairs of sequences $(f_n^-,f_n^+)_{n\in\mathbb{N}}$ such that
$f_n^+\in{\cal C}_n, (f_j^+)_n^+ =f_n^+, (f_j^-)_n^-=f_n^-$ for every $n\in\mathbb{N} \ \ j\geq n$ and $\rho(f_n^-,f_n^+)\rightarrow 0$ as $n\rightarrow\infty.$
\end{Def}

Let
$f\in{\cal C}.$ Then, the pair of sequences $(f_n^-,f_n^+)_{n\in\mathbb{N}}$ belongs to $\widetilde{\cal C}$ (because $f\in{\cal C}_k$ for
some $k\in\mathbb{N}$ and $f_i^-=f_i^+=f$ for $i\geq k).$ This defines an embedding of ${\cal C}$ into $\widetilde{\cal C}$ ; we denote its
image by the same symbol ${\cal C}.$
\\

\begin{Def} We define the following {\sl metric} on the set $\widetilde{\cal C}$ :
$$\widetilde\rho(f,g) = \underset{n\in\mathbb{N}}{\rm Sup} \max(\rho(f_n^-,g_n^-),\rho(f_n^+,g_n^+)).$$
\end{Def}

\begin{Th} $(\widetilde{\cal C},\widetilde\rho)$ {\sl is a complete metric space and ${\cal C}$ is a dense subset of
$\widetilde{\cal C}.$
}
\end{Th}

{\sl Proof :} 
\begin{enumerate}[i)]
\item The completness of $\widetilde{\cal C}.$ Let $(f_{k)}=(f_{(k),n}^-,f_{(k),n}^+)_{n\in\mathbb{N}})_{k\in\mathbb{N}}$ be a Cauchy sequence in
$\widetilde{\cal C}.$ For every fixed $n,$ the sequences $(f_{(k),n}^\pm)_{k\in\mathbb{N}}$ are Cauchy sequences in ${\cal C}_n$ and, by
virtue of the metric completness of ${\cal C}_n,$ converge in ${\cal C}_n$ to limits $F_n^{(-)}$ and $F_n^{(+)}$ respectively.
Moreover, for any $\varepsilon>0,$ there exists $N(\varepsilon)\in\mathbb{N}$ (independent of  $n)$ such that for every $k>N(\varepsilon)$
there are the inequalities
$$\rho(f_{(k),n}^-,F_n^{(-)})<\varepsilon \hbox{ and } \rho(f_{(k),n}^+,F_n^{(+)})<\varepsilon.$$
It remains to show that the pair of sequences. $F = (F_n^{(-)},F_n^{(+)})_{n\in\mathbb{N}}$ belongs to $\widetilde{\cal C}$ ; then the latter
inequalities will imply that $F = \displaystyle\lim_{l\rightarrow\infty}f_{(k)}$ in $\widetilde{\cal C}.$ The equalities
$(F_n^{(\pm)})_j^\pm = F_j^{(\pm)} (n>j)$ follow from the continuity of the maps $f\rightarrow f_j^\pm.$ Finally,
$$\rho(F_n^{(-)},F_n^{(+)})\leq\rho(f_{(k),n}^-,F_n^{(-)})+\rho(f_{(k),n}^+,F_n^{(+)})+\rho(f_{(k),n}^-,f_{(k),n}^+).$$
Hence for $k>N(\varepsilon/3)$ and for sufficiently large $n$ the right-hand side of the latter inequality is less then $\varepsilon.$

\item  ${\cal C}$ is dense in $\widetilde{\cal C}.$ Let $f = (f_n^-,f_n^+)_{n\in\mathbb{N}}\in\widetilde{\cal C}.$ We set
$$F_k = ((f_k^-)_n^-, (f_k^-)_n^+)_{n\in\mathbb{N}}.$$
Each $F_k$ is the image of $f_k^-\in{\cal C}$ under the embedding of ${\cal C}$ into $\widetilde{\cal C}.$ Let us show that the
sequence of $F_k$ converges ot $f$ in the metric $\widetilde\rho.$ We have
\begin{eqnarray*}
&&\widetilde\rho(F_k,f) {\rm Sup}\{\rho(f_1^-,f_1^-),\ldots,\rho(f_k^-,f_k^-), \\
&&\rho(f_k^-,f_{k+1}^-), \rho(f_k^-,f_{k+2}^-),\ldots,\rho((f_k^-)_1^+,f_1^+),\\
&&\rho((f_k^-)_2^+,f_2^+),\ldots,\rho(f_k^-,f_k^+),\rho(f_k^-,f_{k+1}^+),\ldots\}.
\end{eqnarray*}

Let $\varepsilon>0.$ Since $(f_k)_{k\in\mathbb{N}}$ is a Cauchy sequence, all terms of the form $\rho(f_k^-,f_{k+\ell}^+)$ do not exceed
$\varepsilon$ for sufficiently large $k.$ According to condition ii) and the definition of $\widetilde{\cal C},$ all terms of the
form $\rho(f_k^-,f_{k+\ell}^+)$ do not exceed $\varepsilon$ for sufficiently large $k$. It remains to estimate the terms of
the double sequence
$$w_{s,k} = \rho((f_k^-)_s^+,f_s^+),$$
where $k\in\mathbb{N}, s=1,\ldots,k.$
\end{enumerate}

We have to show that there exists $N\in\mathbb{N}$ such that
$$\rho((f_k^-)_s^+,f_s^+)<\varepsilon \eqno(A1)$$
for any $k>N$ and $s\leq k.$ The monotonicity of the mapping $f\rightarrow f_s^+$ and the inequality $f_k^-\leq f_{k+1}^-\leq f_s^+$ show
that $w_{s,k}\geq w_{s,k+1}$ for all $s$ and $k\geq s.$ In particular , $w_{s,k}\leq w_{s,s}$ fr all $k\geq s.$ The sequence
$w_{s,s}$ tends to zero as $s\rightarrow\infty.$ Therefore it is sufficient to proove the inequality (A1) for $k>N$
 and $s\leq s_0,$ where $s_0$ is some fixed number.

The projections $f\rightarrow f_s^+$ are supposed uniformly continuous, then $\exists\delta>0$ such that
$w_{s,k}=\rho((f_k^-)_s^+,(f_k^+)_s^+)<\varepsilon$ if $\rho(f_k^-,f_k^+)<\delta$ and $s\leq s_0, s\leq k.$ Hence the inequality
$w_{s,k}<\varepsilon$ holds for sufficiently large $k$ and (A1) is true. \hfill$\square$\\

\end{document}